\documentclass[12pt]{article}
\usepackage{mathrsfs}
\usepackage{amsmath,amsfonts,amssymb,rotating}
\usepackage{hyperref}

\def\qed{\nopagebreak\hfill{\rule{4pt}{7pt}}}
\def\proof{\noindent {\it{Proof.} \hskip 2pt}}

\newtheorem{theo}{Theorem}[section]

\newtheorem{lemm}[theo]{Lemma}

\newtheorem{coro}[theo]{Corollary}

\numberwithin{equation}{section}

\newdimen\Squaresize \Squaresize=11pt
\newdimen\Thickness \Thickness=0.7pt
\def\Square#1{\hbox{\vrule width \Thickness
   \vbox to \Squaresize{\hrule height \Thickness\vss
    \hbox to \Squaresize{\hss#1\hss}
   \vss\hrule height\Thickness}
\unskip\vrule width \Thickness} \kern-\Thickness}

\def\Vsquare#1{\vbox{\Square{$#1$}}\kern-\Thickness}

\def\moins{\raise 1pt\hbox{{$\scriptstyle -$}}}

\begin{document}

\begin{center}
{\large\bf Infinitely Log-monotonic
Combinatorial Sequences}
\end{center}

\begin{center}
William Y. C. Chen$^{1}$, Jeremy J. F. Guo$^{2}$ and Larry X. W. Wang$^{3}$ \\[8pt]
$^{1,2}$Center for Applied Mathematics\\
Tianjin University\\
 Tianjin 300072, P. R. China\\[6pt]

$^{1,3}$Center for Combinatorics, LPMC-TJKLC\\
Nankai University\\
 Tianjin 300071, P. R. China\\[6pt]

Email: $^{1}${\tt chenyc@tju.edu.cn}, $^{2}${\tt guo@tju.edu.cn}, $^{3}${\tt
wsw82@nankai.edu.cn}
\end{center}
\vspace{0.3cm} \noindent{\bf Abstract.}  We introduce the notion of infinitely  log-monotonic sequences. By establishing a connection between completely monotonic functions and
infinitely log-monotonic sequences, we show that the sequences of the Bernoulli numbers,  the Catalan numbers and the central binomial coefficients are infinitely log-monotonic.  In particular, if a sequence $\{a_n\}_{n\geq 0}$ is log-monotonic of order two, then it is ratio log-concave in the sense that the sequence $\{a_{n+1}/a_{n}\}_{n\geq 0}$
is log-concave.  Furthermore, we prove that if a sequence  $\{a_n\}_{n\geq k}$ is
ratio log-concave, then the sequence
$\{\sqrt[n]{a_n}\}_{n\geq k}$ is strictly log-concave subject to a certain initial condition. As consequences,
we show that the sequences of the derangement numbers, the Motzkin numbers, the Fine numbers, the central Delannoy numbers, the numbers of tree-like polyhexes and the Domb numbers are ratio log-concave. For the case of the Domb numbers $D_n$,
we confirm a conjecture of Sun on the log-concavity of the sequence $\{\sqrt[n]{D_n}\}_{n\geq 1}$.




\section{Introduction}

In this paper, we introduce the notion of infinitely log-monotonic sequences based on the classical concept of
logarithmically completely monotonic functions.
A function $f$ is said to be {\it completely monotonic} on an interval $I$ if $f$ has derivatives of all orders on $I$ and
\begin{equation}
  (-1)^nf^{(n)}(x)\geq 0
\end{equation}
for $x\in I$ and all  integers $n\geq 0$.

A positive function $f$ is said to be {\it logarithmically completely monotonic} on an interval $I$ if $\log f$ satisfies
\begin{equation}
  (-1)^n[\log f(x)]^{(n)}\geq 0
\end{equation}
for $x\in I$ and all the integers $n\geq 1$, see, for example, Atanassov and Tsoukrovski \cite{ata}. A logarithmically completely monotonic function is completely monotonic, but the converse is not necessarily the case, see Berg \cite{ber}.


Recall that a sequence $\{a_n\}_{n\geq 0}$ is said to be {\it log-concave} (resp. {\it log-convex}) if for all $n\geq 1$, $a_n^2\geq a_{n-1}a_{n+1}$ (resp. $a_n^2\leq a_{n-1}a_{n+1}$), and it is said to be {\it strictly} log-concave (resp. {\it strictly} log-convex) if the inequality is strict. Define an operator $\mathcal{R}$ on a sequence $\{a_n\}_{n\geq 0}$ by
\[
\mathcal{R}\{a_n\}_{n\geq 0}=\{b_n\}_{n\geq 0},
\]
where $b_n=a_{n+1}/a_n$. We say that a sequence $\{a_n\}_{n\geq 0}$ is log-monotonic of order $k$ if for $r$  odd and not greater than $k-1$, the sequence $\mathcal{R}^r\{a_n\}_{n\geq 0}$ is log-concave and for $r$   even and
not greater than $k-1$,  the sequence $\mathcal{R}^r\{a_n\}_{n\geq 0}$ is log-convex. A sequence $\{a_n\}_{n\geq 0}$ is called {\it infinitely log-monotonic} if it is log-monotonic of order $k$ for all   integers $k\geq 1$.

We establish a connection between completely monotonic functions and infinitely log-monotonic sequences. Using the log-behavior of the Riemann zeta function and the gamma function, we show that the sequences of the Bernoulli numbers, the Catalan numbers and the central binomial coefficients are infinitely log-monotonic.

Log-monotonic sequences of order two are of special interest. It can be
easily seen that a sequence $\{a_n\}_{n\geq 0}$ is log-monotonic of order two if and only if it is log-convex and  the ratio sequence $\{a_{n+1}/a_n\}_{n\geq 0}$ is log-concave.  A sequence $\{a_n\}_{n\geq 0}$ is said to be ratio log-concave if $\{a_{n+1}/a_n\}_{n\geq 0}$ is log-concave. Similarly, a sequence $\{a_n\}_{n\geq 0}$  is called ratio log-convex if the ratio sequence $\{a_{n+1}/a_n\}_{n\geq 0}$ is log-convex.

We prove that under a certain initial condition, the
ratio log-concavity of a sequence $\{a_{n}\}_{n\geq k}$ of positive numbers implies
that the sequence $\{ \sqrt[n]{a_n}\}_{n \geq k}$ is strictly log-concave.
Analogous to Firoozbakht's conjecture that the sequence $\{\sqrt[n]{p_n}\}_{n\geq 1}$ is
strictly increasing, where $p_n$ is the $n$-th prime number,
Sun \cite{sun} conjectured that for some combinatorial sequences $\{a_n\}_{n\geq 0}$, the
sequences $\{\sqrt[n]{a_n}\}_{n\geq 1}$ are strictly log-concave except for the first few terms. Hou, Sun and Wen \cite{hou}
proved this conjecture for the derangement numbers. Luca and St\u{a}nic\u{a}'s \cite{luc} showed that the conjecture holds for the Bernoulli numbers, the Euler numbers, the Tangent numbers, the Motzkin numbers, the Ap\'{e}ry numbers, the Franel numbers, the central Delannoy numbers, the Schr\"{o}der numbers and the trinomial coefficients.

In this paper, we prove that for
the derangement numbers, the Motzkin numbers, the Fine numbers, the central Delannoy numbers, the numbers of tree-like polyhexes and the Domb numbers, the sequences $\{a_n\}_{n\geq 0}$ are ratio
log-concave except for the first few terms. As
consequences, we are led to the log-concavity of the sequence $\{\sqrt[n]{a_n}\}_{n\geq 1}$ for the derangement numbers proved by Hou, Sun and Wen and the log-concavity of the sequences $\{\sqrt[n]{a_n}\}_{n\geq 1}$ for the Bernoulli numbers, the Motzkin numbers and the central Delonnay numbers proved by  Luca and St\u{a}nic\u{a}. Moreover,
we confirm the conjecture of Sun on the log-concavity of the sequence $\{\sqrt[n]{D_n}\}_{n\geq 1}$ of the Domb numbers. We also show that for the
sequences $\{a_n\}_{n\geq 0}$ of the Catalan numbers, the central binomial coefficients, the Fine numbers and the numbers of tree-like polyhexes, the sequences $\{\sqrt[n]{a_n}\}_{n\geq 0}$ are strictly log-concave except for the first few terms.

Our approach also applies to the log-behavior of the harmonic numbers. We prove that for any given integer $m\geq 1$, the
sequence of the generalized harmonic numbers $\{H_{n,m}\}_{n\geq 1}$ is ratio
log-convex. This leads to the strict log-convexity of the sequence $\{\sqrt[n]{H_{n,m}}\}_{n\geq 3}$, as conjectured by Sun \cite{sun} and confirmed by Hou, Sun and Wen \cite{hou}.

We conclude this paper with a conjecture
on almost infinitely log-monotonic sequences,
where a sequence $\{a_n\}_{n\geq 0}$ is called almost infinitely log-monotonic if for $k\geq 0$, $\{a_n\}_{n\geq 0}$ is log-monotonic
of order $k$ except for certain entries at the beginning.

\section{Infinitely Log-monotonic Sequences}

In this section, we establish a connection between completely
 monotonic functions and  infinitely log-monotonic sequences.
 We show that
 the Riemann zeta function is logarithmically completely monotonic.
 Then we deduce that the sequence of the Bernoulli numbers is
 infinitely log-monotonic.
Using the fact that $[\log \Gamma(x)]''$ is completely monotonic, we prove that the sequences of the Catalan numbers and the central binomial coefficients are infinitely log-monotonic.

The following theorem shows that a completely monotonic function
gives rise to an infinitely log-monotonic sequence.

\begin{theo}\label{rlcm}
Assume that  $f(x)$ is a function such that $[\log f(x)]''$ is completely monotonic for $x\geq 1$. Let $a_n=f(n)$ for $n\geq 1$. Then  the sequence $\{a_n\}_{n\geq 1}$ is infinitely log-monotonic.
\end{theo}

\proof Denote the sequence $\mathcal{R}^{i}\{a_n\}_{n\geq 1}$ by $\{b_{n,i}\}_{n\geq 1}$ for  $i\geq 0$. In other words,
$b_{n,0}=a_n$  and $b_{n,i+1}=b_{n+1,i}/b_{n,i}$.
Let $f_0(x)=f(x)$, and define the functions $f_1(x),\ f_2(x),\ldots$ by
the relation
\begin{equation}\label{lem11}
f_{i+1}(x)={ f_i(x+1)\over f_i(x)}.
\end{equation}
It is clear that  $b_{n,i}=f_i(n)$ for any $i\geq 0$ and $n\geq 1$.
We aim to show that for $j\geq 0$ and $k\geq 2$,
\begin{equation}\label{111}
(-1)^k[\log f_{2j}(x)]^{(k)}\geq 0,
\end{equation}
and
\begin{equation}\label{222}
(-1)^k[\log f_{2j+1}(x)]^{(k)}\leq 0.
\end{equation}

We use induction on $j$. Since $[\log f(x)]''$ is  completely monotonic, we
find that for $k\geq 2$,
\begin{equation}\label{2}
(-1)^k[\log f(x)]^{(k)}\geq 0,
\end{equation}
that is, \eqref{111} holds for $j=0$.
Rewriting \eqref{2} as $(-1)^{k+1}[\log f(x)]^{(k+1)}\geq 0$ for $k\geq 1$,
we get
\begin{equation}\label{lem22}(-1)^{k}[\log f(x)]^{(k+1)}\leq 0,
\end{equation}
for $k\geq 1$. Since $f_1(x)=f(x+1)/f(x)$, by \eqref{lem22} we find that
for $k\geq 2$,
\[
(-1)^k[\log f_1(x)]^{(k)}=(-1)^k[\log f(x+1)]^{(k)}-(-1)^k[\log f(x)]^{(k)}\leq 0.
\]
Thus  \eqref{222} is verified for $j=0$.

We now assume that \eqref{111} and \eqref{222} hold for  $j\leq n-1$.
We proceed to show that \eqref{111} and \eqref{222} hold for $j=n$.
Rewriting the induction hypothesis \eqref{222} as   $(-1)^{k+1}[\log f_{2n-1}(x)]^{(k+1)}\leq 0$ for $k\geq 1$, we see  that for $k\geq 2$,
\[
 (-1)^k[\log f_{2n}(x)]^{(k)}=(-1)^k[\log f_{2n-1}(x+1)]^{(k)}-(-1)^k[\log f_{2n-1}(x)]^{(k)}\geq 0.
\]
Hence \eqref{111} holds for $j=n$. Similarly, it can be shown that
for $k\geq 2$,
 \[ (-1)^k[\log f_{2n+1}(x)]^{(k)}\leq 0,\]
  that is, \eqref{222} holds for $j=n$. Up to now, we have proved \eqref{111} and \eqref{222} for $j\geq 0$ and $k\geq 2$.

In view of \eqref{111} and \eqref{222}, we conclude that for any $i\geq 0$, the sequence $\{f_{2i}(n)\}_{n\geq 1}$ is log-convex and the sequence $\{f_{2i+1}(n)\}_{n\geq 1}$ is log-concave, in other words,
 for any $i\geq 0$, the sequence $\mathcal{R}^{2i}\{a_n\}_{n\geq 1}$ is log-convex and the sequence $\mathcal{R}^{2i+1}\{a_n\}_{n\geq 1}$ is log-concave. This completes the proof.\qed


It is not difficult to see that the  Riemann zeta function
is logarithmically completely monotonic.
Indeed, it is known that
for $\Re(s)>1$,
\begin{equation}\label{lzeta}
-\frac{\zeta^{\prime}(s)}{\zeta(s)}=\sum_{n=1}^{\infty}\Lambda(n)n^{-s},
\end{equation}
where $\Lambda(n)$ is the von Mangoldt function defined by
\begin{equation}\label{mangoldt}
 \Lambda(n)=\left\{
\begin{array}{ll}
     \log p, &\mbox{if }n=p^m,\mbox{ where $p$ is a prime and $m\geq 1$,}\\[3pt]
     0,      &\mbox{otherwise,}
\end{array}
\right.
\end{equation}
see \cite[pp.3]{ivi}.
Using the above formula \eqref{lzeta}, we find  that, for $x>1$,
\[
 (-1)^k[\log \zeta(x)]^{(k)}
 =\sum_{n=1}^{\infty}\frac{\Lambda(n)(\log n)^{k-1}}{n^x}>0.
\]
It follows that $\zeta(x)$ is logarithmically completely monotonic on $(1,\infty)$. Thus $[\log \zeta(x)]''$ is completely monotonic on $(1,\infty)$.

To prove that the sequence of  the Bernoulli numbers $\{|B_{2n}|\}_{n\geq 1}$ is
infinitely log-monotonic, we need to show that $[\log \Gamma(x)]''$ is completely monotonic on $(0,\infty)$. It is known that for $x>0$,
\begin{equation}
 \frac{\Gamma^{\prime}(x)}{\Gamma(x)}= -\gamma-\sum_{m=1}^{\infty}\left(\frac{1}{x+m-1}-\frac{1}{m}\right),
\end{equation}
 where $\gamma$ is the Euler constant, see  \cite{alz}. It follows that
\begin{equation}\label{psi}
  \left( \frac{\Gamma^{\prime}(x)}{\Gamma(x)}\right)^{(k)}=(-1)^{k+1}k!\sum_{m=0}^{\infty}\frac{1}{(x+m)^{k+1}}.
\end{equation}
Therefore, for $k\geq 2$ and $x\geq 1$, we have $(-1)^k[\log \Gamma(x)]^{(k)}>0$, namely, $[\log \Gamma(x)]''$ is completely monotonic on $(0,\infty)$.

\begin{theo} Let $B_n$ be the $n$-th Bernoulli number.
  The sequence $\{|B_{2n}|\}_{n\geq 1}$ is infinitely log-monotonic.
\end{theo}

\proof 
It is well known  that
\begin{equation}\label{rmzf}
  \zeta(2n)=\frac{2^{2n-1}\pi^{2n}}{(2n)!}|B_{2n}|.
\end{equation}
 Set
\begin{equation}\label{theta}
  z(x) =\frac{2\zeta(2x)\Gamma(2x+1)}{(2\pi)^{2x}}.
\end{equation}
Using  \eqref{rmzf}, it can be checked that $z(n)=|B_{2n}|$.
Since for $x\geq 1$ and $k\geq 2$,
\[
[\log z(x)]^{(k)}=2^k[\log \zeta(2x)]^{(k)}+2^k[\log \Gamma(2x+1)]^{(k)},
\]
we see  that for $x\geq 1$ and  $k\geq 2$, $(-1)^k[\log z(x)]^{(k)}>0$. Thus for $x\geq 1$, $[\log z(x)]''$ is completely monotonic. Hence the proof is complete by using Theorem \ref{rlcm}.
\qed

Next we consider the Catalan numbers $C_n$.
Define \[c(x)=\frac{\Gamma(2x+1)}{\Gamma(x+1)\Gamma(x+2)},\] so that $c(n)=C_n=\frac{1}{n+1}{{2n}\choose n}$ for $n\geq 1$. To prove that the sequence $\{C_n\}_{n\geq 1}$ is infinitely log-monotonic, we shall
 demonstrate that $[\log c(x)]''$ is completely monotonic on $[1,\infty]$.

\begin{theo}\label{catalan}
The sequence $\{C_n\}_{n\geq 1}$ of the Catalan numbers is infinitely log-monotonic.
\end{theo}
\proof By the definition of $c(x)$, we find that for $x\geq 1$ and $k\geq 2$,
\begin{equation}\label{ca}
[\log c(x)]^{(k)}=[\log \Gamma(2x+1)]^{(k)}-[\log\Gamma(x+1)]^{(k)}-[\log\Gamma(x+2)]^{(k)}.
\end{equation}
In view of \eqref{psi}, we obtain that for $x\geq 1$ and $k\geq 2$,
\begin{eqnarray}
    \lefteqn{(-1)^k [\log \Gamma(2x+1)]^{(k)}}\nonumber\\
 &=&2^k(k-1)!\sum_{m=0}^{\infty}\frac{1}{(2x+1+m)^k}\nonumber\\
 &=&(k-1)!\sum_{i=0}^{\infty}\frac{1}{(x+1/2+(2i)/2)^k}+(k-1)!\sum_{i=0}^{\infty}\frac{1}{(x+1/2+(2i+1)/2)^k}\nonumber\\
 &=&(k-1)!\sum_{i=0}^{\infty}\frac{1}{(x+1/2+i)^k}+(k-1)!\sum_{i=0}^{\infty}\frac{1}{(x+i+1)^k}.\label{ca2}
\end{eqnarray}
Since  \[(k-1)!\sum_{i=0}^{\infty}\frac{1}{(x+1/2+i)^k}>(k-1)!\sum_{m=0}^{\infty}\frac{1}{(x+1+m)^k}
=(-1)^k[\log\Gamma(x+1)]^{(k)}\]
and
\[(k-1)!\sum_{i=0}^{\infty}\frac{1}{(x+i+1)^k}>(k-1)!\sum_{m=0}^{\infty}\frac{1}{(x+2+m)^k}
=(-1)^k[\log\Gamma(x+2)]^{(k)},\]
it follows from \eqref{ca2} that
\[
{(-1)^k [\log \Gamma(2x+1)]^{(k)}}>(-1)^k[\log\Gamma(x+1)]^{(k)}+(-1)^k[\log\Gamma(x+2)]^{(k)}.
\]
According to \eqref{ca}, we find that for $x\geq 1$ and $k\geq 2$, $(-1)^k[\log c(x)]^{(k)}>0$, that is, $[\log c(x)]''$ is completely monotonic for $x\geq 1$.
By Theorem \ref{rlcm}, we conclude that the sequence $\{C_n\}_{n\geq 1}$ is infinitely log-monotonic.\qed

For the sequence  $\{{{2n}\choose n}\}_{n\geq 1}$ of the central binomial coefficients, we define
\[
d(x)=\frac{\Gamma(2x+1)}{\Gamma(x+1)^2}.
\]
 Then we have $d(n)={{2n}\choose n}$ for $n\geq 1$. Using the same argument
 as in the proof of Theorem \ref{catalan}, it can be shown that for $x\geq 1$, $(-1)^k[\log d(x)]^{(k)}>0$. By Theorem \ref{rlcm}, we come to the conclusion that the sequence $\{{{2n}\choose n}\}_{n\geq 1}$ is infinitely log-monotonic.

\section{Ratio Log-Concavity}

In this section, we  show that the ratio log-concavity (resp. ratio log-convexity) of a sequence $\{a_n\}_{n\geq k}$ implies the strict log-concavity (resp. strict log-convexity) of the sequence $\{\sqrt[n]{a_n}\}_{n\geq k}$ under a certain initial condition.
We also show that the sequence of the derangement numbers is ratio log-concave and the sequence of the generalized harmonic numbers is ratio log-convex. Some known results are
consequences of the log-behavior of these two ratio sequences.

\begin{theo}\label{rlcav}
Assume that $k$ is a positive integer.
  If a sequence $\{a_n\}_{n\geq k}$ is ratio log-concave and
  \begin{equation}\label{kk3}
  \frac{\sqrt[k+1]{a_{k+1}}}{\sqrt[k]{a_k}}> \frac{\sqrt[k+2]{a_{k+2}}}{\sqrt[k+1]{a_{k+1}}},\end{equation}  then the sequence $\{\sqrt[n]{a_n}\}_{n\geq k}$ is strictly log-concave.
\end{theo}

To prove above theorem, we need the following lemmas.

\begin{lemm}\label{rll}
Assume that $k$ is a positive integer. If a sequence $\{a_n\}_{n\geq k}$ is ratio log-concave, then for any $k\leq i< j \leq n$, we have
\begin{equation}\label{ineq1}
    \log\frac{a_{j+1}}{a_j}
    \geq (j-i)\left(\log \frac{a_{n+1}}{a_n}-\log \frac{a_n}{a_{n-1}}\right)+\log\frac{a_{i+1}}{a_i}.
\end{equation}
\end{lemm}

\proof The ratio log-concavity of $\{a_n\}_{n\geq k}$ implies that
for  $k\leq i\leq n-1$,
\[ \left(\frac{a_{i+2}}{a_{i+1}}\right)^2\geq \frac{a_{i+1}}{a_i}\cdot \frac{a_{i+3}}{a_{i+2}}.\]
Thus, for $k\leq i\leq n-1$,
\[\log \frac{a_{i+2}}{a_{i+1}}-\log\frac{a_{i+1}}{a_{i}} \geq
    \log\frac{a_{i+3}}{a_{i+2}}-\log\frac{a_{i+2}}{a_{i+1}}.
\]
It follows that for $k\leq i\leq n-1$,
\begin{equation}\label{kk7}
   \log \frac{a_{i+2}}{a_{i+1}}-\log\frac{a_{i+1}}{a_{i}} \geq
    \log\frac{a_{n+1}}{a_{n}}-\log\frac{a_{n}}{a_{n-1}}.
\end{equation}
Since for $j>i$,
\[
    \log\frac{a_{j+1}}{a_j}
    =\sum_{l=i}^{j-1}\left(\log\frac{a_{l+2}}{a_{l+1}}-\log\frac{a_{l+1}}{a_{l}}\right)+\log\frac{a_{i+1}}{a_i},
\]
by \eqref{kk7}, we finally  get \eqref{ineq1}.
\qed

Under the initial condition \eqref{kk3} in Theorem \ref{rlcav}, we have the following inequality.

\begin{lemm}\label{rll2}
Assume that $k$ is a positive integer. If a ratio log-concave sequence $\{a_n\}_{n\geq k}$ satisfies the initial
condition \eqref{kk3} in Theorem \ref{rlcav},
then we have for $n>k$,
\begin{equation}\label{22}
 k\log\frac{a_{k+1}}{a_k}>    \frac{k^2+k}{2}\left(\log\frac{a_{n+1}}{a_n}-\log \frac{a_n}{a_{n-1}}\right)+\log a_k.
\end{equation}
\end{lemm}

\proof The initial condition \eqref{kk3} can be rewritten as
\[
    \frac{\log a_{k+1}}{k+1}-\frac{\log a_k}{k}>\frac{\log a_{k+2}}{k+2}-\frac{\log a_{k+1}}{k+1},
\]
that is,
\begin{equation} \label{kl}
    2(k^2+2k)\log a_{k+1}>(k^2+k)\log a_{k+2}+(k^2+3k+2)\log a_k.
\end{equation}
Expressing (\ref{kl}) in terms of logarithms of ratios,
we find that
\begin{equation}\label{kk10}
    k\log\frac{a_{k+1}}{a_k}> \frac{k^2+k}{2}\left(\log\frac{a_{k+2}}{a_{k+1}}-\log\frac{a_{k+1}}{a_k}\right)+\log a_k.
\end{equation}
By  the ratio log-concavity of $\{a_n\}_{n\geq k}$, we deduce that
for $k<n$,
\begin{equation}\label{kl2}
\log\frac{a_{k+2}}{a_{k+1}}-\log\frac{a_{k+1}}{a_k}\geq 
\cdots \geq \log\frac{a_{n+1}}{a_n}-\log\frac{a_n}{a_{n-1}}.
\end{equation}
Combining \eqref{kk10} and \eqref{kl2}, we arrive at \eqref{22}.
\qed

We are now ready to complete the proof of Theorem \ref{rlcav}.

\noindent {\it  Proof of Theorem \ref{rlcav}.} To establish the strict log-concavity of
$\{\sqrt[n]{a_n}\}_{n\geq k}$, we aim to show that for $n> k$,
\begin{equation}
\frac{2\log a_n}{n}> \frac{\log a_{n+1}}{n+1}+\frac{\log a_{n-1}}{n-1}.
\end{equation}
It is easily verified that
\begin{eqnarray}
  \lefteqn{ \frac{2\log a_n}{n}-\frac{\log a_{n+1}}{n+1}-\frac{\log a_{n-1}}{n-1}}\nonumber\\
  &=& \frac{1}{n(n+1)(n-1)}(2(n^2-1)\log a_n-(n^2-n)\log a_{n+1}-(n^2+n)\log a_{n-1})\nonumber\\[5pt]
  &=& \frac{1}{n^3-n}\left[n\left(\log \frac{a_{n+1}}{a_n}+\log \frac{a_n}{a_{n-1}}\right)-n^2\left(\log \frac{a_{n+1}}{a_n}-\log \frac{a_n}{a_{n-1}}\right)-2\log a_n
  \right].\nonumber\\ \label{main1}
\end{eqnarray}
We wish to estimate $n\log \frac{a_{n+1}}{a_n}$ and $n\log \frac{a_n}{a_{n-1}}$. We claim that
\begin{equation}\label{444}
    n\log\frac{a_{n+1}}{a_n}
    > \frac{n(n+1)}{2}\left(\log \frac{a_{n+1}}{a_n}-\log \frac{a_n}{a_{n-1}}\right)+\log a_n,
\end{equation}
and
\begin{equation}\label{555}
    n\log\frac{a_n}{a_{n-1}}
    > \frac{n(n-1)}{2}\left(\log \frac{a_{n+1}}{a_n}-\log \frac{a_n}{a_{n-1}}\right)
     +\log a_n.
\end{equation}
Setting $j=n$ in \eqref{ineq1} of Lemma \ref{rll}, we see that for $k\leq i\leq n-1$,
\begin{equation}\label{kk13}
    \log\frac{a_{n+1}}{a_n}
    \geq (n-i)\left(\log \frac{a_{n+1}}{a_n}-\log \frac{a_n}{a_{n-1}}\right)+\log\frac{a_{i+1}}{a_i}.
\end{equation}
Summing \eqref{kk13} over $i$ from $k$ to $n-1$ gives
\begin{equation}\label{333}
    (n-k)\log\frac{a_{n+1}}{a_n}\geq \sum_{i=k}^{n-1}(n-i)\left(\log \frac{a_{n+1}}{a_n}-\log \frac{a_n}{a_{n-1}}\right)+\sum_{i=k}^{n-1}\log\frac{a_{i+1}}{a_i}.
\end{equation}
Write $k\log\frac{a_{n+1}}{a_n}$ in the following form
\begin{equation}\label{kk11}
 k\log\frac{a_{n+1}}{a_n}
 =k\sum_{j=k+1}^n\left(\log\frac{a_{j+1}}{a_j}
 -\log\frac{a_j}{a_{j-1}}\right)+k\log\frac{a_{k+1}}{a_k}.
\end{equation}
It follows from the ratio log-concavity condition
 \eqref{kl2} that
\begin{equation}\label{23}
 k\sum_{j=k+1}^n\left(\log\frac{a_{j+1}}{a_j}-\log\frac{a_j}{a_{j-1}}\right)
 \geq  k(n-k)\left(\log \frac{a_{n+1}}{a_n}-\log \frac{a_n}{a_{n-1}}\right).
\end{equation}
Applying \eqref{23} and Lemma \ref{rll2} to the righthand side of \eqref{kk11}, we deduce that
\begin{eqnarray}
  k\log\frac{a_{n+1}}{a_n}&>& \left( k(n-k)+\frac{k^2+k}{2} \right)\left(\log \frac{a_{n+1}}{a_n}-\log \frac{a_n}{a_{n-1}}\right)+\log a_k\nonumber\\
   & & = \sum_{i=0}^{k-1}(n-i)\left(\log \frac{a_{n+1}}{a_n}-\log \frac{a_n}{a_{n-1}}\right)+\log a_k.\label{ineq2}
\end{eqnarray}
Combining \eqref{333} and \eqref{ineq2}, we obtain that
\begin{eqnarray*}
    n\log\frac{a_{n+1}}{a_n}
    &>& \sum_{i=0}^{n-1}(n-i)\left(\log \frac{a_{n+1}}{a_n}-\log \frac{a_n}{a_{n-1}}\right)+\sum_{i=k}^{n-1}\log\frac{a_{i+1}}{a_i}+\log a_k \nonumber\\
    & & = \frac{n(n+1)}{2}\left(\log \frac{a_{n+1}}{a_n}-\log \frac{a_n}{a_{n-1}}\right)+\log a_n.
\end{eqnarray*}
This verifies \eqref{444}.

We continue to prove \eqref{555}.
Setting $j=n-1$ in \eqref{ineq1} in Lemma \ref{rll}, we find  that for $k\leq i\leq n-2$,
\begin{equation}\label{kk12}
    \log\frac{a_n}{a_{n-1}}
    \geq (n-i-1)\left(\log \frac{a_{n+1}}{a_n}-\log \frac{a_n}{a_{n-1}}\right)+\log\frac{a_{i+1}}{a_i}.
\end{equation}
Summing \eqref{kk12} over $i$ from $k$ to $n-2$, we get
\begin{equation}\label{kk4}
(n-k-1)\log\frac{a_n}{a_{n-1}}\geq \sum_{i=k}^{n-2}(n-i-1)\left(\log \frac{a_{n+1}}{a_n}-\log \frac{a_n}{a_{n-1}}\right)+\sum_{i=k}^{n-2}\log\frac{a_{i+1}}{a_i}.
\end{equation}
Note that
\begin{equation}\label{kk15}
k\log\frac{a_n}{a_{n-1}}=k\sum_{j=k+1}^{n-1}\left(\log\frac{a_{j+1}}{a_j}
-\log\frac{a_j}{a_{j-1}}\right)+k\log\frac{a_{k+1}}{a_k}.
\end{equation}
Using the ratio log-concavity condition
 \eqref{kl2}, we find that
\begin{equation}\label{kk8}
 k\sum_{j=k+1}^{n-1}\left(\log\frac{a_{j+1}}{a_j}-\log\frac{a_j}{a_{j-1}}\right)
 \geq  k(n-k-1)\left(\log \frac{a_{n+1}}{a_n}-\log \frac{a_n}{a_{n-1}}\right).
\end{equation}
Applying \eqref{kk8} and Lemma \ref{rll2} to the righthand side of \eqref{kk15}, we
obtain that
\begin{equation}\label{kk5}
    k\log\frac{a_n}{a_{n-1}} > \sum_{i=0}^k(n-i-1)\left(\log \frac{a_{n+1}}{a_n}-\log \frac{a_n}{a_{n-1}}\right)+\log a_k.
\end{equation}
Combining \eqref{kk4} and \eqref{kk5}, we have
\begin{eqnarray*}
    (n-1)\log\frac{a_n}{a_{n-1}}
  &>& \sum_{i=0}^{n-2}(n-i-1)\left(\log \frac{a_{n+1}}{a_n}-\log \frac{a_n}{a_{n-1}}\right)
        +\sum_{i=k}^{n-2}\log\frac{a_{i+1}}{a_i}+\log a_k\\
  & & \quad =\ \frac{n(n-1)}{2}\left(\log \frac{a_{n+1}}{a_n}-\log \frac{a_n}{a_{n-1}}\right)
     +\log a_{n-1}.
\end{eqnarray*}
This proves \eqref{555}.

Utilizing  the estimates of $n\log \frac{a_{n+1}}{a_n}$ and $n\log \frac{a_n}{a_{n-1}}$ as given in \eqref{444} and \eqref{555}, we
deduce that
\[
    n\left(\log \frac{a_{n+1}}{a_n}+\log \frac{a_n}{a_n-1}\right) >
    n^2\left(\log \frac{a_{n+1}}{a_n}-\log \frac{a_n}{a_{n-1}}\right)+2\log a_n.
\]
It follows from \eqref{main1}  that
\[
   \frac{2\log a_n}{n}>\frac{\log a_{n+1}}{n+1}+\frac{\log a_{n-1}}{n-1}.
\]
This completes the proof.\qed

In Section 2, we have shown that the sequences $\{|B_{2n}|\}_{n\geq 1}$,
$\{C_n\}_{n\geq 1}$ and $\{{{2n}\choose n}\}_{n\geq 1}$ are infinitely
log-monotonic, and hence they are ratio log-concave.
For the sequences $\{\sqrt[n]{|B_{2n}|}\}_{n\geq 2}$, $\{\sqrt[n]{C_n}\}_{n\geq 1}$ and $\{{\sqrt[n]{{2n}\choose n}}\}_{n\geq 1}$, it can be checked that
the initial condition \eqref{kk3} in Theorem \ref{rlcav} holds. Thus we obtain the following
properties.

\begin{coro}
The sequences $\{\sqrt[n]{|B_{2n}|}\}_{n\geq 2}$, $\{\sqrt[n]{C_n}\}_{n\geq 1}$ and $\{{\sqrt[n]{{2n}\choose n}}\}_{n\geq 1}$ are strictly log-concave.
\end{coro}


Next we show that the sequence of the derangement numbers are ratio log-concave.
For $n\geq 1$, the $n$-th derangement number $D_n$ is the number of permutations $\sigma$ of $\{1, 2, \ldots, n\}$ that
have no fixed points, that is, $\sigma(i)\neq i$ for  $i=1, 2, \ldots, n$. It is known that the sequence $\{D_n\}_{n\geq 2}$ is log-convex, see Liu and Wang \cite{liu}.

\begin{theo}\label{derangement}
The sequence $\{D_n\}_{n\geq 2}$ is ratio log-concave.
\end{theo}

\proof To prove that $\{D_n\}_{n\geq 2}$ is ratio log-concave, we
 proceed to verify that
\begin{equation}\label{ll3}
\frac{D_{n+2}^2}{D_{n+1}^2}>\frac{D_{n+3}}{D_{n+2}}\cdot \frac{D_{n+1}}{D_n}.
\end{equation}
Using the recurrence relation
\begin{equation}\label{derrec}
  D_n= nD_{n-1}+ (-1)^n,
\end{equation}
we get
\begin{eqnarray}
 \lefteqn {\frac{D_{n+2}^2}{D_{n+1}^2}-\frac{D_{n+3}}{D_{n+2}}\frac{D_{n+1}}{D_n}}\\
 &=&\left(n+2+\frac{(-1)^{n+2}}{D_{n+1}}\right)^2-
 \left(n+3+\frac{(-1)^{n+3}}{D_{n+2}}\right)\left(n+1+\frac{(-1)^{n+1}}{D_n}\right)\nonumber\\
 & & \geq 1-\frac{n+1}{D_{n+2}}-\frac{n+3}{D_n}-\frac{2(n+2)}{D_{n+1}}
 +\frac{1}{D_{n+1}^2}-\frac{1}{D_{n+2}D_n}.\label{dd}
\end{eqnarray}
From \eqref{derrec}, it is easily seen that for $n\geq 5$, \begin{equation}\label{d1}
D_n>5(n+3).
\end{equation}
It follows that
\begin{equation}\label{ll}
\frac{n+1}{D_{n+2}}+\frac{n+3}{D_n}+\frac{2(n+2)}{D_{n+1}}<\frac{4}{5}.
\end{equation}
Since the sequence $\{D_n\}_{n\geq 2}$ is log-convex, that is,
\begin{equation}\label{ll2}
\frac{1}{D_{n+1}^2}>\frac{1}{D_{n+2}D_n},
\end{equation}
applying \eqref{ll}  to the righthand side of \eqref{dd}, we deduce that \eqref{ll3} holds for $n\geq 5$. On the other hand, it is easily verified that \eqref{ll3} holds for $2\leq n\leq 5$.
Hence the proof is complete. \qed

Since $\sqrt[4]{D_4}/\sqrt[3]{D_3}>\sqrt[5]{D_5}/\sqrt[4]{D_4}$, by Theorem \ref{rlcav} we
 are led to the known result that the sequence $\{\sqrt[n]{D_n}\}_{n\geq 3}$
is strictly log-concave.


The above approach also applies to ratio log-convex sequences. We have the following criterion.

\begin{theo}\label{rlvex}
 Assume that $k$ is a positive integer. If a sequence $\{a_n\}_{n\geq k}$ is ratio log-convex and
   \[\frac{\sqrt[k+1]{a_{k+1}}}{\sqrt[k]{a_k}}< \frac{\sqrt[k+2]{a_{k+2}}}{\sqrt[k+1]{a_{k+1}}},\] then the sequence $\{\sqrt[n]{a_n}\}_{n\geq k}$ is strictly log-convex.
\end{theo}

As an application of the above theorem, we find that for any $m\geq 1$,  the sequence
 $\{H_{m,n}\}_{n\geq 1}$ of the generalized harmonic numbers is ratio log-convex.
Recall that for any positive integers $m,\ n$, the $n$-th generalized harmonic number $H_{n,m}$ of order $m$ is defined by
\begin{equation}\label{hl}
 H_{n,m}=\sum_{k=1}^n \frac{1}{k^m}.
\end{equation}

\begin{theo}\label{har}
 For any $m\geq 1$, the sequence $\{H_{n,m}\}_{n\geq 1}$ is ratio log-convex.
\end{theo}

\proof Towards the assertion of the theorem, we aim to show that $\frac{H_{n+2,m}H_{n,m}}{H_{n+1,m}^2}$ is strictly increasing in $n$.
By the definition given by \eqref{hl}, we get the following recurrence relations
\begin{equation}\label{h1}
H_{n+2,m}=H_{n+1,m}+\frac{1}{(n+2)^k}
\end{equation}
and
\begin{equation}\label{h2}
H_{n,m}=H_{n+1,m}-\frac{1}{(n+1)^k}.
\end{equation}
Thus
\begin{eqnarray}
\frac{H_{n+2,m}H_{n,m}}{H_{n+1,m}^2}
&=&1-\left(\frac{1}{(n+1)^m}-\frac{1}{(n+2)^m}\right)\frac{1}{H_{n+1,m}}\nonumber\\
& &\quad-\frac{1}{(n+1)^m(n+2)^mH_{n+1,m}^2}.\label{h3}
\end{eqnarray}
Since for $x>0$, $(x^{-m}-(x+1)^{-m})$ is strictly decreasing, we see that $\frac{1}{n^m}-\frac{1}{(n+1)^m}$ is strictly decreasing in $n$.
By \eqref{h3}, we deduce that $\frac{H_{n+2,m}H_{n,m}}{H_{n+1,m}^2}$ is strictly increasing in $n$, that is, the sequence $\{H_{n,m}\}_{n\geq 1}$ is ratio log-convex for any $m\geq 1$.\qed

For any $m\geq 1$, it can be checked that $\sqrt[4]{H_{4,m}}/\sqrt[3]{H_{3,m}}<\sqrt[5]{H_{5,m}}/\sqrt[4]{H_{4,m}}$.
As a consequence of Theorem \ref{rlvex} and Theorem \ref{har}, we arrive at the known result that the sequence $\{\sqrt[n]{H_{n,m}}\}_{n\geq 3}$ is strictly log-convex for any $m\geq 1$.

\section{Sequences Satisfying a Three-term Recurrence Relation}

In this section, we are concerned with sequences $\{a_n\}_{n\geq 0}$
satisfying a three-term recurrence relation
\begin{equation}\label{3term}
  a_n=u(n)a_{n-1}+v(n)a_{n-2}
\end{equation}
where $u(n)$ and $v(n)$ are  rational functions in $n$ and $u(n)>0$ for $n\geq 2$. We shall give a criterion  for $\{a_n\}_{n\geq k}$
to be ratio log-concave.  Using this criterion, it
can be shown that the sequences of the Motzkin numbers, the Fine numbers, the central Delannoy numbers, the numbers of the tree-like polyhexes and the Domb numbers are ratio log-concave. As consequences, we are led to two known results for the Motzkin numbers and the central Delannoy numbers.
Moreover, we  confirm a conjecture of Sun \cite{sun} on the log-behavior of the Domb numbers.

First, we consider the case when $v(n)>0$ for any $n\geq 2$.

\begin{theo}\label{3termrc}
Let $\{a_n\}_{n\geq 0}$ be the sequence defined by the recurrence relation \eqref{3term}. Assume that for $n\geq 2$, $v(n)>0$  and $u(n)^3>u(n+1)v(n)$. If there exists a nonnegative integer $N$ and a function $g(n)$ such that for all $n\geq N+2$,
\begin{itemize}
\item[(i)] $a_n/a_{n-1}\geq g(n) \geq u(n)$;
\item[(ii)] $g(n)^4-u(n) g(n)^3-u(n+1)v(n) g(n)-v(n)v(n+1)> 0$,
\end{itemize}
then the sequence $\{a_n\}_{n\geq N}$ is ratio log-concave.

\end{theo}

\proof To prove that the sequence $\{a_n\}_{n\geq N}$ is ratio log-concave, we proceed to show that for $n\geq N+2$,
\begin{equation}\label{41}
a_n^3a_{n-2}-a_{n+1}a_{n-1}^3>0.
\end{equation}
By the recurrence relation \eqref{3term}, we deduce that
\begin{eqnarray*}
  \lefteqn{a_n^3a_{n-2}-a_{n+1}a_{n-1}^3}\\
    &=&\frac{1}{v(n)}a_n^3(a_n-u(n) a_{n-1})-(u(n+1)a_n+v(n+1)a_{n-1})a_{n-1}^3\\
    &=&\frac{a_{n-1}^4}{v(n)}\left[\left(\frac{a_n}{a_{n-1}}\right)^4-u(n)\left(\frac{a_n}{a_{n-1}}\right)^3\right.\\
    & & \ \left.-u(n+1)v(n)\left(\frac{a_n}{a_{n-1}}\right)-v(n)v(n+1)\right].
\end{eqnarray*}
Since $v(n)>0$ for $n\geq 2$, in order to prove \eqref{41}, it suffices to verify that for $n\geq N+2$,
\begin{equation}\label{lt}
  \left(\frac{a_n}{a_{n-1}}\right)^4-u(n)\left(\frac{a_n}{a_{n-1}}\right)^3
        -u(n+1)v(n)\left(\frac{a_n}{a_{n-1}}\right)-v(n)v(n+1)>0.
\end{equation}
Define  \[f(x)=x^4-u(n)x^3-u(n+1)v(n)x-v(n)v(n+1).\]
Then (\ref{lt}) can be rewritten as  $f(\frac{a_n}{a_{n-1}})>0$.
Note that
\[
 f'(x)=4x^3-3u(n)x^2-u(n+1)v(n)
\]
and
\[
 f''(x)=12x^2-6u(n)x.
\]
Hence $f''(x)>0$ for $x>u(n)/2$. This implies that for $x>u(n)/2$, $f'(x)$ is strictly increasing. Using the condition  $u(n)^3>u(n+1)v(n)$ in the theorem, we get $f'(u(n))>0$. It follows that for $x\geq u(n)$, $f'(x)>0$. Thus $f(x)$ is increasing for $x\geq u(n)$. Given the condition $g(n)\geq u(n)$, we see that  $f(x)$ is strictly increasing for $x\geq g(n)$. On the other hand, condition (ii)  says that $f(g(n))> 0$ for any $n\geq N+2$.
So we have $f(x)>0$ for $x\geq g(n)$. From the condition  $a_n/a_{n-1}\geq g(n)$ in the theorem, we deduce that $f(\frac{a_n}{a_{n-1}})>0$ for $n\geq N+2$. This completes the proof.\qed

To apply the above theorem, we need a lower bound $g(n)$ on the
ratio $a_{n}/a_{n-1}$ subject to conditions (i) and (ii).
The following lemma will be used to give a heuristic approach
to finding a lower bound $g(n)$.

\begin{lemm}\label{gn1}
Let $\{a_n\}_{n\geq 0}$ be the sequence defined by the recurrence relation \eqref{3term}. Assume that $v(n)>0$ for $n\geq 2$. If there exists a positive integer $N$ and a function $g(n)$ such that
\[ g(N)<\frac{a_{N}}{a_{N-1}}<\frac{v(N+1)}{g(N+1)-u(N+1)}\] and the inequalities
\begin{equation}\label{gn}
u(n)+\frac{v(n)}{g(n-1)}<\frac{v(n+1)}{g(n+1)-u(n+1)}
\end{equation}
hold for all $n\geq N$, then for $n\geq N$,
\begin{equation}\label{l2}
g(n)<\frac{a_n}{a_{n-1}}<\frac{v(n+1)}{g(n+1)-u(n+1)}.
\end{equation}
\end{lemm}

\proof We use induction on $n$. Assume that \eqref{l2} hold for $n$, where $n\geq N$.
We proceed to show that \eqref{l2} also holds for $n+1$, that is,
\begin{equation}\label{lll4}
g(n+1)<\frac{a_{n+1}}{a_{n}}<\frac{v(n+2)}{g(n+2)-u(n+2)}.
\end{equation}
By recurrence relation \eqref{3term}, we find that
\[
\frac{a_{n+1}}{a_n}=\frac{u(n+1)a_n+v(n+1)a_{n-1}}{a_n}=u(n+1)+v(n+1)\frac{a_{n-1}}{a_n}.
\]
So \eqref{lll4} is equivalent to
\begin{equation}\label{lll5}
g(n+1)<u(n+1)+v(n+1)\frac{a_{n-1}}{a_n}<\frac{v(n+2)}{g(n+2)-u(n+2)}.
\end{equation}
But the second inequality
 \[\frac{a_n}{a_{n-1}}<\frac{v(n+1)}{g(n+1)-u(n+1)}\]
  in the induction hypothesis can be rewritten as
\begin{equation}\label{lll2}
 u(n+1)+v(n+1)\frac{a_{n-1}}{a_n}>g(n+1),
\end{equation}
this yields    the first inequality in \eqref{lll5}.
On the other hand, using condition \eqref{gn} with $n$ replaced by
 $n+1$, we get
\begin{equation}\label{lll6}
u(n+1)+\frac{v(n+1)}{g(n)}<\frac{v(n+2)}{g(n+2)-u(n+2)},
\end{equation}
By the first inequality $g(n)<a_{n}/a_{n-1}$ in the induction hypothesis and  the above inequality \eqref{lll6}, we find that
\begin{equation}\label{lll3}
u(n+1)+v(n+1)\frac{a_{n-1}}{a_n}
<u(n+1)+\frac{v(n+1)}{g(n)}<\frac{v(n+2)}{g(n+2)-u(n+2)}.
\end{equation}
This leads to the second inequality in \eqref{lll5}. Hence the proof is complete by induction.\qed

Notice that Theorem \ref{3termrc} and Lemma \ref{gn1} require  a lower bound $g(n)$ of $a_n/a_{n-1}$ subject to certain conditions. Adopting the framework in \cite{chen1}, we
give an iterative procedure to find a lower bound $g(n)$. It should be noted that
the success of this procedure is not guaranteed, although it serves the
purpose in many cases.

Let us begin with the quadratic equation
\begin{equation}\label{lambda}
 \lambda^2-u(n)\lambda-v(n)=0.
\end{equation}
Since $v(n)>0$ for $n\geq 2$, this equation has a unique positive root
\begin{equation}\label{M2}
 \lambda(n)=\frac{u(n)+\sqrt{u(n)^2+4v(n)}}{2}.
\end{equation}
If $\lambda(n)$ satisfies condition (ii) in Theorem \ref{3termrc} and \eqref{gn}, then it is a feasible choice for $g(n)$. Otherwise, we
construct a function $r(n)$ such that $r(n)>\lambda(n)$ for any nonnegative integer $n$. Since $u(n)$ and $v(n)$ are rational functions, we may assume that
\begin{equation}\label{mm1}
u(n)^2+4v(n)=\frac{P(n)}{Q(n)},
\end{equation}
where $P(n)$ and $Q(n)$ are polynomials in $n$. If $P(n)$ can be written as $R(n)^2-c$, where
$R(n)$ is a polynomial in $n$ and $c$ is a positive number, then we have
\[\sqrt{u(n)^2+4v(n)}<\frac{R(n)}{\sqrt{Q(n)}}.\]
Clearly,  $r(n)$ can be
chosen as follows to meet the requirement $r(n)>\lambda(n)$,
\begin{equation}\label{sn}
r(n)=\frac{u(n)\sqrt{Q(n)}+R(n)}{2\sqrt{Q(n)}}.
\end{equation}
If $r(n)$ satisfies condition (ii) in Theorem \ref{3termrc} and \eqref{gn}, then it is the desired lower bound. Otherwise, we try to find a
number $x$ such that
\begin{equation}\label{f11}
g(n)=r(n)+\frac{1}{d(n)}\frac{x}{n}
\end{equation}
satisfies condition (ii) in Theorem \ref{3termrc} and \eqref{gn}, where $d(n)$ is the denominator of $r(n)$.
Since the lower bound $g(n)$ in Lemma \ref{gn1} satisfies the
 two inequalities in  \eqref{l2}, this implies that
\begin{equation}\label{pp1}
\frac{v(n+1)}{g(n+1)-u(n+1)}>g(n).
\end{equation}
We shall be guided by the above inequality \eqref{pp1} in search for the number $x$. More precisely, let
\begin{equation}\label{cxn}
  C(x,n)=\frac{v(n+1)}{g(n+1)-u(n+1)}-g(n),
\end{equation}
where $g(n)$ is given by \eqref{f11}, and
 \[g(n+1)=r(n+1)+\frac{1}{d(n+1)}\frac{x}{n+1}.\]
If $C(x,n)$ is a rational function, then let
\[
C(x,n)=\frac{Y(x,n)}{Z(x,n)},
\]
where $Y(x,n)$ and $Z(x,n)$ are polynomials in $x$ and $n$.
 We now consider $x$ as a number and treat $Y(x,n)$ as a polynomial $Y(n)$ in $n$. Denote by $H(x)$ the coefficient of the term of   highest degree in $Y(n)$, and set $H(x)=0$. If $x_1$ is a solution of $H(x)=0$, then we set \[g(n)=r(n)+\frac{1}{d(n)}\frac{x_1}{n}.\] If $g(n)$ satisfies condition (ii) in Theorem \ref{3termrc} and \eqref{gn}, then it is the desired lower bound.   Otherwise, we repeat the above process to find a number $x_2$ such that \[g(n)=r(n)+\frac{1}{d(n)}\left(\frac{x_1}{n}+\frac{x_2}{n^2}\right)\] satisfies condition (ii) in Theorem \ref{3termrc} and \eqref{gn}. If we are lucky, by iteration we may find numbers $\ x_1,\ x_2,\ \ldots,\ x_k$ such that \[g(n)=r(n)+\frac{1}{d(n)}\left(\frac{x_1}{n}+\frac{x_2}{n^2}+\cdots+\frac{x_k}{n^k}\right)\] satisfies the lower bound condition (ii) in Theorem \ref{3termrc} and \eqref{gn}.
This leads to a lower bound $g(n)$ of the ratio $a_{n}/a_{n-1}$.

For example, let us consider the Motzkin numbers $M_n$ defined by the recurrence relation
\begin{equation}\label{M}
M_n=\frac{2n+1}{n+2}M_{n-1} +\frac{3n-3}{n+2}M_{n-2},
\end{equation}
where $n\geq 2$ and  $M_0=M_1=1$, see Aigner \cite{aig}. In the context of the general recurrence relation \eqref{3term}, for the case of $M_n$, we have $u(n)=(2n+1)/(n+2)$ and $v(n)=(3n-3)/(n+2)$. It is easy to see that  for $n\geq 2$, $v(n)>0$ and $u(n)^3>u(n+1)v(n)$, that is, $v(n)$ and $u(n)$ satisfy the conditions in Theorem \ref{3termrc}.

Following the
procedure of finding a lower bound $g(n)$ of $a_n/a_{n-1}$, we begin with
 the unique positive root of equation \eqref{lambda}
\[
\lambda(n)=\frac{2n+1+\sqrt{16n^2+16n-23}}{2(n+2)}.
\]
Clearly, $\lambda(n)$ does not satisfy condition (ii) in Theorem \ref{3termrc}. So we continue to construct a function $r(n)$ such that $r(n)> \lambda(n)$ for any positive integer $n$.
Since
\[
u(n)^2+4v(n)=\frac{(4n+2)^2-27}{(n+2)^2},
\]
we have $R(n)=4n+2$ and $Q(n)=(n+2)^2$. By \eqref{sn}, $r(n)$ can be chosen as
\[
r(n)=\frac{3n+\frac{3}{2}}{n+2}.
\]
It can be checked that $r(n)$ does not satisfy inequality \eqref{gn} in Lemma \ref{gn1}. Then we try to find a number $x$ such that
\[g(n)=\frac{3n+\frac{3}{2}}{n+2}+\frac{x}{n+2}\]
satisfies condition (ii) in Theorem \ref{3termrc} and \eqref{gn}. By the definition \eqref{cxn} of $C(x,n)$, we have
\[C(x,n)=\frac{-(16x+9)n^2-(16x+9)n-4x^2-6x}{2(2n^2+5n+3+2x)(n+2)n}.\] Let $Y(x,n)$ denote the numerator of $C(x,n)$. Setting the coefficient of $n^2$ in $Y(x,n)$ to zero,
we get $-16x-9=0$ and  $x_1=-\frac{9}{16}$.
So we get \[g(n)=\frac{6n^2+3n-\frac{9}{8}}{2n(n+2)}.\]
It is easy to verify that for $n\geq 13$, $g(n)$ satisfies the conditions in Theorem \ref{3termrc} and Lemma \ref{gn1}. Thus we  deduce that $\{M_n\}_{n\geq 11}$ is ratio log-concave. Moreover, it can be verified that for $6\leq n\leq 12$, \[M_n^3M_{n-2}>M_{n+1}M_{n-1}^3.\]
Hence we arrive at the following assertion.

\begin{theo}\label{mot}
The sequence $\{M_n\}_{n\geq 4}$ of the Motzkin numbers  is ratio log-concave.
\end{theo}

Since for $2\leq n\leq 5$, $(\sqrt[n]{M_n})^2>\sqrt[n+1]{M_{n+1}}\sqrt[n-1]{M_{n-1}}$,
  Theorem \ref{rlcav} and Theorem \ref{mot} imply the known result that the sequence $\{\sqrt[n]{M_n}\}_{n\geq 1}$ is strictly log-concave.

The Fine numbers $f_n$ are given by the recurrence relation
\begin{equation}\label{fine}
2(n + 1)f_n = (7n-5)f_{n-1} + 2(2n-1)f_{n-2},
\end{equation}
where $n\geq 2$ and $f_0=1$ and $f_1=0$, see Deutsch and Shapiro \cite{deu}. Next we show that the sequence $\{f_n\}_{n\geq 5}$ is ratio log-concave.

\begin{theo}
  The sequence $\{f_n\}_{n\geq 5}$ of the Fine numbers is ratio log-concave, and the sequence $\{\sqrt[n]{f_n}\}_{n\geq 2}$ is strictly log-concave.
\end{theo}

\proof For the Fine numbers $f_n$, we have  $u(n)=(7n-5)/(2n+2)$ and $v(n)=(2n-1)/(n+1)$ in the recurrence relation  \eqref{3term}. Clearly, for $n\geq 2$ we have $v(n)>0$ and $u(n)^3>u(n+1)v(n)$, that is, $v(n)$ and $u(n)$ satisfy the conditions in Theorem \ref{3termrc} with $n\geq 2$. Employing the above procedure to find
a lower bound $g(n)$ of $f_n/f_{n-1}$, we get
 \[g(n)=\frac{4n^2-2n+\frac{2}{3}}{n^2+n}.\]
It can be  verified that for $n\geq 7$, $g(n)$ satisfies conditions in Theorem \ref{3termrc} and Lemma \ref{gn1}. By Theorem \ref{3termrc}, we deduce that $\{f_n\}_{n\geq 5}$ is ratio log-concave.
Moreover, for $2\leq n\leq 5$,  we have \[(\sqrt[n]{f_n})^2>\sqrt[n+1]{f_{n+1}}\sqrt[n-1]{f_{n-1}}.\] By Theorem \ref{rlcav}, we
conclude that $\{\sqrt[n]{f_n}\}_{n\geq 1}$ is strictly log-concave.
 \qed

Parallel to Theorem \ref{3termrc}, we give a criterion for the ratio log-concavity
of a sequence satisfying a three-term recurrence relation
\begin{equation}\label{3tt}
  a_n=u(n)a_{n-1}+v(n)a_{n-2}
\end{equation}
where $v(n)<0$ for $n\geq 2$. The proof of the following theorem is
analogous to the proof of Theorem 4.1, and hence is omitted.

\begin{theo}\label{3termrc'}
Let $\{a_n\}_{n\geq 0}$ be the sequence defined by the recurrence relation \eqref{3tt}. Assume that $v(n)<0$ for $n\geq 2$. If there exists a nonnegative
integer $N$ and a function $h(n)$ such that for all $n\geq N+2$,
\begin{itemize}
\item[(i)] $3u(n)/4\leq a_n/a_{n-1}\leq h(n)$;
\item[(ii)] $h(n)^4-u(n)h(n)^3-u(n+1)v(n)h(n)-v(n)v(n+1)< 0$,
\end{itemize}
then $\{a_n\}_{n\geq N}$ is ratio log-concave.
\end{theo}

The following lemma will be used to derive
an upper  bound $h(n)$ of the ratio $a_n/a_{n-1}$, which is
needed in Theorem \ref{3termrc'}. The proof of this lemma is omitted since it is
essentially the same as the proof of Lemma \ref{gn1}.

\begin{lemm}\label{hn1}
Let $\{a_n\}_{n\geq 0}$ be the sequence defined by the recurrence relation \eqref{3tt}. Assume that $v(n)<0$ for $n\geq 2$. If there exists a positive integer $N$ and a function $h(n)$ such that $a_N/a_{N-1}<h(N)$ and the equalities
  \begin{equation}\label{hn}
    h(n+1)>u(n+1)+\frac{v(n+1)}{h(n)},
  \end{equation}
  hold for all $n\geq N$, then for $n\geq N$, we have $a_n/a_{n-1}<h(n)$.
\end{lemm}

We now describe an iterative procedure  to find an upper bound $h(n)$
of $a_n/a_{n-1}$ which is a rational
function in $n$ satisfying certain conditions. Since $v(n)<0$ for any $n\geq 2$, equation \eqref{lambda} has  either no real roots
 or two positive roots. If
  \eqref{lambda} has two positive roots, then $\lambda(n)$ given by \eqref{M2} is the larger root. If $\lambda(n)$ satisfies condition (ii) in Theorem \ref{3termrc'} and \eqref{hn}, then it is a feasible choice for $h(n)$. Otherwise, we
construct a function $s(n)$ such that $s(n)<\lambda(n)$ for any nonnegative integer $n$. If $P(n)$ given by \eqref{mm1} can be written as $T(n)^2+c$, where
$T(n)$ is a polynomial in $n$ and $c$ is a positive number, then
we may choose $s(n)$ to be
\begin{equation}\label{rn}
s(n)=\frac{u(n)\sqrt{Q(n)}+T(n)}{2\sqrt{Q(n)}},
\end{equation}
where $Q(n)$ is given by \eqref{mm1}.
If $s(n)$ satisfies  condition (ii) in Theorem \ref{3termrc'} and \eqref{hn}, then it is a feasible upper bound. Otherwise, we try to find a
number $x$ such that
\begin{equation}\label{f22}
h(n)=s(n)+\frac{1}{d(n)}\frac{x}{n}
\end{equation}
satisfies condition (ii) in Theorem \ref{3termrc'} and \eqref{hn}, where $d(n)$ is the denominator of $s(n)$.
We use inequality \eqref{hn} to look for the number $x$. Let
\begin{equation}\label{dxn}
  D(x,n)=h(n+1)-u(n+1)-\frac{v(n+1)}{h(n)},
\end{equation}
where $h(n)$ is given by \eqref{f22}, and
 \[h(n+1)=r(n+1)+\frac{1}{d(n+1)}\frac{x}{n+1}.\]
If $D(x,n)$ is a rational function, then we obtain $x_1$ from $D(x,n)$
in the same manner as we get $x_1$ from $C(x,n)$ in the above
procedure of deriving $g(n)$ as a lower bound of $a_n/a_{n-1}$. We set \[h(n)=s(n)+\frac{1}{d(n)}\frac{x_1}{n}.\] If $h(n)$ satisfies condition (ii) in Theorem \ref{3termrc'} and \eqref{hn}, then it is a feasible upper bound. Otherwise, we repeat the above procedure  to find a number $x_2$ such that \[h(n)=s(n)+\frac{1}{d(n)}\left(\frac{x_1}{n}+\frac{x_2}{n^2}\right)\] satisfies  condition (ii) in Theorem \ref{3termrc'} and \eqref{hn}. Eventually, by using this process we may find numbers $\ x_1,\ x_2,\ \ldots,\ x_k$ such that \[h(n)=s(n)+\frac{1}{d(n)}\left(\frac{x_1}{n}+\frac{x_2}{n^2}+\cdots+\frac{x_k}{n^k}\right)\] satisfies the upper bound condition (ii) in Theorem \ref{3termrc'} and \eqref{hn}. Then we get a required bound $h(n)$ of the ratio $a_{n}/a_{n-1}$.


For example,   consider the central Delannoy numbers $D(n)$ defined by the recurrence relation
 \begin{equation}
  D(n)=\frac{3(2n-1)}{n}D(n-1)-\frac{n-1}{n}D(n-2),
\end{equation}
where $n\geq 2$ and $D(0)=1$ and $D(1)=3$, see Sun \cite{sun1}.

In the context of   recurrence relation \eqref{3tt}, for the case of $D(n)$, we have $u(n)=3(2n-1)/n$ and $v(n)=-(n-1)/n$. It can be seen that for $n\geq 2$, $v(n)<0$ and \[ {D(n)\over D(n-1)} \geq {3u(n)\over 4},
    \]
 that is, $u(n)$ and $v(n)$ meet the requirements
 of Theorem \ref{3termrc'} with $n\geq 2$. In this case, the quadratic equation \eqref{lambda} becomes
 \begin{equation}\label{cdn}
 \lambda^2-\frac{3(2n-1)}{n}\lambda+\frac{n-1}{n}=0,
 \end{equation}
Thus the larger root of equation \eqref{cdn} is
 \[
 \lambda(n)=\frac{6n-3+\sqrt{32n^2-32n+9}}{2n}.
 \]
It can be seen that $\lambda(n)$ does not satisfy condition (ii) in Theorem \ref{3termrc'}. Now, we wish to find a function $s(n)$ such that $s(n)<\lambda(n)$ for any positive integer $n$. In this case,
\[
u(n)^2+4v(n)=\frac{2(4n-2)^2+1}{n^2},
\]
so we may set $T(n)=4\sqrt{2}n-2\sqrt{2}$ and $Q(n)=n^2$. By \eqref{rn}, $s(n)$ can be chosen as
\[
s(n)=\frac{(3+2\sqrt{2})n-\frac{3}{2}-\sqrt{2}}{n}.
\]
It can be checked that $s(n)$ does not satisfy condition (ii) in Theorem \ref{3termrc'}. So we further consider \[h(n)=\frac{(3+2\sqrt{2})n-\frac{3}{2}-\sqrt{2}}{n}+\frac{1}{n}\frac{x}{n}.\]
By the definition \eqref{dxn} of $D(x,n)$, we have
\[
  D(x,n)=\frac{(16\sqrt{2}x+1)n^2-(24x-8\sqrt{2}x-1)n+4\sqrt{2}x-6x+4x^2}
  {2(n+1)(6n^2+4\sqrt{2}n^2-2\sqrt{2}n-3n+2x)}.
\]
Let $Y(x,n)$ be the numerator of $D(x,n)$, which is a
polynomial of degree $2$ in $n$. Setting the coefficient of $n^2$ in $Y(x,n)$ to zero, we obtain $16\sqrt{2}x+1=0$, $x_1=-\sqrt{2}/32$ and
\[h(n)=\frac{(3+2\sqrt{2})n^2-\frac{3}{2}n-\sqrt{2}n-\frac{\sqrt{2}}{32}}{n^2}.\]
At this time,  $h(n)$ satisfies the conditions in Theorem \ref{3termrc'} and Lemma \ref{hn1}.  By Theorem \ref{3termrc'}, we reach the following conclusion.

\begin{theo}\label{denn}
The sequence $\{D(n)\}_{n\geq 0}$ of the central Delannoy numbers is ratio log-concave.
\end{theo}

Since  $D(2)>D(1)\sqrt[3]{D(3)}$,    Theorem \ref{rlcav} and Theorem \ref{denn} lead to the known result that the sequence $\{\sqrt[n]{D(n)}\}_{n\geq 1}$ is strictly log-concave.

Let $t_n$ be the number of  tree-like polyhexes with $n+1$ hexagons, which is given by the recurrence relation
\[
 (n+1)t_n=3(2n-1)t_{n-1}-5(n-2)t_{n-2},
\]
where $n\geq 2$ and $t_0=t_1=1$, see Harary and Read \cite{hara}. The
following theorem shows that the sequence $\{t_n\}_{n\geq 0}$
is ratio log-concave.

\begin{theo}
The sequence $\{t_n\}_{n\geq 2}$ is ratio log-concave, and the sequence $\{\sqrt[n]{t_n}\}_{n\geq 1}$ is strictly log-concave.
\end{theo}

\proof For the numbers of tree-like polyhexes, we have $u(n)=(6n-3)/(n+1)$ and $v(n)=-(5n-10)/(n+1)$ in the recurrence \eqref{3tt}. Clearly, for $n\geq 2$ we have $v(n)<0$ and
$t_n/t_{n-1}\geq 3u(n)/4$,
that is, $u(n)$ and $v(n)$ satisfy the conditions in Theorem \ref{3termrc'} with $n\geq 2$. Using the above procedure to find
an upper bound $h(n)$ of the ratio $t_n/t_{n-1}$, we get \[h(n)=\frac{10n^3-5n^2+\frac{15}{8}n+6}{2n^2+2n^3}.\] For $n\geq 7$, it can be verified that $h(n)$ satisfies the conditions in Theorem \ref{3termrc'} and Lemma \ref{hn1}. Thus by Theorem \ref{3termrc'}, we deduce that $\{t_n\}_{n\geq 5}$ is ratio log-concave. Moreover, for $2\leq n\leq 6$, we have \[t_n^3t_{n-2}-t_{n+1}t_{n-1}^3>0.\] Thus the sequence $\{t_n\}_{n\geq 0}$ is ratio log-concave.

 It can be verified that $\sqrt{t_2}/t_1>\sqrt[3]{t_3}/\sqrt{t_2}$.
 Hence  Theorem \ref{rlcav} implies that the sequence $\{\sqrt[n]{t_n}\}_{n\geq 1}$ is strictly log-concave.  \qed

Finally, we consider the sequence of the Domb numbers $D_n$ given by the recurrence relation
\[
 n^3D_n=2(2n-1)(5n^2-5n+2)D_{n-1}-64(n-1)^3D_{n-2},
\]
where $n\geq 2$ and $D_0=1$ and $D_1=4$. The $n$-th Domb number $D_n$
is the number of $2n$-step polygons on the diamond lattice. Chan,  Chan and Liu \cite{cha} obtained a series for $\frac{1}{\pi}$ involving the Domb numbers. Chan, Tanigawa, Yangc, and Zudilin \cite{cha2} found three analogues of Clausen's identities involving Domb numbers.
Sun \cite{sun} conjectured that the sequence $\{D_n\}_{n\geq 0}$ is log-convex and the sequence $\{\sqrt[n]{D_n}\}_{n\geq 1}$ is strictly increasing and strictly log-concave. Wang and Zhu \cite{wang} proved the sequence $\{D_n\}_{n\geq 0}$ is log-convex and the sequence $\{\sqrt[n]{D_n}\}_{n\geq 1}$ is increasing.
Next we show $\{D_n\}_{n\geq 0}$ is ratio log-concave. As a consequence, we confirm the conjecture of Sun on the strict log-concavity of $\{\sqrt[n]{D_n}\}_{n\geq 1}$.

We shall slightly modify the above procedure to
derive an upper bound $h(n)$ of $D_n/D_{n-1}$. It turns
out that with some adjustment we can find an upper bound $h(n)$
without iteration.

\begin{theo}
The sequence $\{D_n\}_{n\geq 0}$ of the Domb numbers is ratio log-concave, and the sequence $\{\sqrt[n]{D_n}\}_{n\geq 1}$ is strictly log-concave.
\end{theo}

\proof As far as the  recurrence \eqref{3tt} is concerned,
for the Domb numbers, we have $u(n)=2(2n-1)(5n^2-5n+2)/n^3$ and $v(n)=-64(n-1)^3/n^3$. It is easy to verify that for $n\geq 2$, $v_n<0$ and for $n\geq 24$, $D_n/D_{n-1}>3u(n)/4$, that is, $u(n)$ and $v(n)$ satisfy the conditions in Theorem \ref{3termrc'} with $n\geq 24$.
By the definition \eqref{M2} of $\lambda(n)$, we have
\begin{eqnarray*}
\lambda(n)&=&\frac{20n^3-30n^2+18n-4}{2n^3}\\
& &\quad+\frac{\sqrt{(12n^3-18n^2+22n-8)^2+208n^2-208n+48}}{2n^3}.
\end{eqnarray*}
It can be verified that $\lambda(n)$ does not satisfy condition (ii) in Theorem \ref{3termrc'}. In order to find a rational function $s(n)$
such that $s(n)<\lambda(n)$ for $n\geq 0$, we may ignore
 $208n^2-208n+48$ in the square root. Set
\[
s(n)=\frac{16n^3-24n^2+40n-12}{n^3}.
\]
Notice that $s(n)$ does not satisfy condition (ii) in Theorem \ref{3termrc'}.
By adjusting the constant term and the coefficient of $n$ in $s(n)$, we get a rational function
\[
h(n)=\frac{16n^3-24n^2+12n-2}{n^3}.
\]
Clearly, $h(n)<s(n)$ for $n\geq 1$. It can be checked that for $n\geq 24$, $h(n)$ satisfies the conditions in Theorem \ref{3termrc'} and Lemma \ref{hn1}. Thus by Theorem \ref{3termrc'} the sequence $\{D_n\}_{n\geq 22}$ is ratio log-concave. Moreover, for $2\leq n\leq 23$ $D_n^3D_{n-2}>D_{n-1}^3D_{n+1}$. Hence, the sequence $\{D_n\}_{n\geq 0}$ is ratio log-concave.

It is easily seen that $\sqrt{D_2}/D_1>\sqrt[3]{D_3}/\sqrt{D_2}$. By Theorem \ref{rlcav},
we are led to the assertion that the sequence $\{\sqrt[n]{D_n}\}_{n\geq 1}$ is strictly log-concave. \qed

To conclude, we conjecture that the sequences of the Motzkin numbers, the Fine numbers, the central Delannoy numbers, the numbers of tree-like polyhexes and the Domb numbers are almost infinitely log-monotonic.
More precisely, we say that a sequence  is
almost infinitely log-monotonic if for each $k\geq 0$,
it is log-monotonic of order $k$ except for certain terms
at the beginning.

We also conjecture that the sequence of the Bell numbers $B_n$ is almost infinitely log-monotonic, where $B_n$ is the number of
partitions of $\{1, 2, \ldots , n\}$.

\vspace{0.5cm}
\noindent{\bf Acknowledgments.}  This work was supported by  the 973
Project, the PCSIRT Project of the Ministry of Education, the Doctoral Program Fund of the Ministry of Education, and the
National Science Foundation of China.

\end{document}